\documentclass[journal,twoside,web]{ieeecolor}
\usepackage{lcsys}

\usepackage{cite}
\usepackage{amsmath,amssymb,amsfonts}
\usepackage{algorithmic}
\usepackage{graphicx}
\usepackage{textcomp}

\usepackage{mathtools}
\usepackage{booktabs}
\usepackage{units}
\usepackage[dvipsnames]{xcolor}
\usepackage{comment}
\usepackage{caption}
\usepackage{hyperref}
\hypersetup{
    colorlinks=true,
    linkcolor=blue,     
    urlcolor=blue,
}
\usepackage{xparse}
\usepackage{orcidlink}

\newcommand{\R}{\mathbb{R}}
\newcommand{\N}{\mathbb{N}}

\newcommand{\mbb}[1]{\mathbb{#1}}
\newcommand{\mrm}[1]{\mathrm{#1}}
\newcommand{\mc}[1]{\mathcal{#1}}
\newcommand{\diag}[1]{\mrm{diag}\big({#1}\big)}

\newcounter{thmcounter}[section]  % Resets with each new section
\renewcommand{\thethmcounter}{\thesection.\arabic{thmcounter}}

\newcommand{\defthmwithqed}[2]{%
  \NewDocumentEnvironment{#1}{ o }{%
    \refstepcounter{thmcounter}% increment counter and allow referencing
    \begin{trivlist}%
      \item[\hskip \labelsep \bfseries #2~\thethmcounter%
        \IfValueT{##1}{\ (##1)}.]%
  }{%
    \hfill $\square$%
    \end{trivlist}%
  }%
}

\newtheorem{theorem}{Theorem}
\newtheorem{lemma}{Lemma}
\newtheorem{proposition}{Proposition}
\newtheorem{corollary}{Corollary}
\newtheorem{defi}{Definition}
\newtheorem{assumption}{Assumption}

\newcommand{\rev}[1]{#1}

\usepackage{tikz}
\usepackage{pgfplots}
\usepackage{hyperref}
\usepgfplotslibrary{groupplots} % for grouping plots
\usetikzlibrary{calc} % for calculations within TikZ
\pgfplotsset{compat=1.18}
\usetikzlibrary{decorations.pathreplacing,decorations.markings}
\usetikzlibrary{positioning}
\usetikzlibrary{backgrounds,scopes}
\usetikzlibrary {arrows.meta}
\usepackage{tikz-network}
\SetLayerDistance{-0.6}
\def\BibTeX{{\rm B\kern-.05em{\sc i\kern-.025em b}\kern-.08em
    T\kern-.1667em\lower.7ex\hbox{E}\kern-.125emX}}
\begin{document}
\bstctlcite{IEEEexample:BSTcontrol}
\title{Exact Time-Varying Turnpikes for \\Dynamic Operation of District Heating Networks}

\author{Max Rose$^{1}$\orcidlink{0000-0002-3985-1326},
Hannes Gernandt$^{1,2}$\orcidlink{0000-0001-7364-4606},
Timm Faulwasser$^{3}$\orcidlink{0000-0002-6892-7406},
Johannes Schiffer$^{1,4}$\orcidlink{0000-0001-5639-4326}
\thanks{$^{1}$Fraunhofer Research Institution for Energy Infrastructures and Geothermal Technologies IEG, 03046 Cottbus, Germany, {\tt\small max.rose@ieg.fraunhofer.de}}%
\thanks{$^{2}$Institute of Mathematical Modelling, Analysis and Computational Mathematics IMACM, University of Wuppertal, 42119 Wuppertal, Germany, {\tt\small gernandt@uni-wuppertal.de}}
\thanks{$^{3}$Institute of Control Systems, Hamburg University of Technology, 21079 Hamburg, Germany, {\tt\small timm.faulwasser@ieee.org}}%
\thanks{$^{4}$Brandenburg University of Technology Cottbus-Senftenberg BTU, 03046 Cottbus, Germany, {\tt\small schiffer@b-tu.de}}
\thanks{MR, HG, and JS acknowledge funding from the German Federal Government, the Federal Ministry of Education and Research, and the State of Brandenburg within the framework of the joint project EIZ: Energy Innovation Center (project numbers 85056897 and 03SF0693A). TF acknowledges support by the Deutsche Forschungsgemeinschaft (DFG, German Research Foundation) (project number 519323897). HG acknowledges support  by the Deutsche Forschungsgemeinschaft (DFG, German
Research Foundation) – Project-ID 531152215 – CRC 1701. 
}
}
\maketitle
\thispagestyle{empty}
\pagestyle{empty}

\begin{abstract}
District heating networks (DHNs) are crucial for decarbonizing the heating sector.
Yet, their efficient and reliable operation requires the coordination of multiple heat producers and the consideration of future demands.
Predictive and optimization-based control is commonly used to address this task, but existing  results for DHNs do not account for  time-varying problem aspects.
Since the turnpike phenomenon can serve as a basis for model predictive control design and analysis, this paper examines its role in DHN optimization by analyzing the underlying optimal control problem with time-varying prices and demands.
That is, we derive conditions for the existence of a unique time-varying singular arc, which constitutes the time varying turnpike, and we provide its closed-form expression.
Additionally, we present converse turnpike results  showing a exact time-varying case implies strict dissipativity of the optimal control problem.
A numerical example illustrates our findings.
\end{abstract}


%--------------------------------------------------------------------------------------------------------------
\section{Introduction}
%--------------------------------------------------------------------------------------------------------------
\IEEEPARstart{I}{n} urban areas of the northern hemisphere, district heating networks (DHNs) are a key technology to tackle climate change through the defossilized heat energy supply~\cite{HeatRoadmap}.
A main driver for this is the potential of DHNs to integrate large amounts of renewable energy to satisfy heat demands of entire cities with minimal CO$_2$-emissions \cite{HeatRoadmap}.
A DHN is essentially a system of pipes used to transport heated water as an energy carrier between producers and consumers. %\cite{HandbookDHN}.

Conventionally, DHNs are operated with just a few, often fossil fuel-based heat producers injecting a controllable heat flow. 
The operating strategies for such DHNs are well understood and are mostly based on simple feedback controls without active coordination of the injected heat flows or online optimization of operating costs~\cite{Buffa2021}. 
In contrast, \rev{in so called $4^\mrm{th}$ generation DHNs}, integrating multiple renewable \rev{energy-based heat} producers requires complex coordination of the injected heat flows to satisfy demands along with physical constraints while achieving an economically profitable operation.
To address this challenge, different optimization-based operating strategies, mainly variants of model predictive control (MPC), are proposed \cite{Buffa2021}.

In \cite{BENONYSSON1995297} and \cite{Rose_IFAC23} an optimization-based economic operation of DHNs is considered.
The formulation of numerically efficient optimal control problems (OCPs) is investigated in \cite{krug_nonlinear_2021} and \cite{rein_optimal_2020}.
The solution of an OCP subject to DHN dynamics modeled as differential-algebraic equations (DAEs) is analyzed in \cite{Jäkle04112024}.
All of the aforementioned works have the common focus of open-loop optimal control. 
However, no MPC-like feedback is considered, and the closed-loop dynamics of the DHN, which results from the application of the receding-horizon solution of the respective OCP to the DHN dynamics, is not investigated.
In contrast, in \rev{\cite{JANSEN2024122874} and \cite{HERING2021120140}, MPC utilizing mixed-integer optimization is considered but no formal guarantees are developed.
MPC with stability guarantees for the closed-loop dynamics is provided in} \cite{Rose_ECC24} and \cite{Sibeijn_ECC24} but the considered scenarios are limited to stationary operating points.
%In applications, 
Yet, due to changing heat flow demands or varying energy prices, it is not to be expected that stationary DHN operation is economically beneficial in real-world applications. 

Put differently, there exists a gap in the available literature concerning economically motivated controls with time-varying demands and prices together with closed-loop guarantees for DHNs in non-stationary settings.
This gap is addressed in the present paper. Our approach is guided by the observation that in \cite{Sibeijn_ECC24} the turnpike phenomenon occurring in an OCP with time-varying costs solved for DHNs is reported, which in turn can be exploited in the closed-loop MPC analysis~\cite{kit:faulwasser18c}.
Hence, by combining these two facts, we conduct a turnpike analysis of the underlying OCP with time-varying cost functions and heat demands for DHNs.

The turnpike phenomenon arises in parametric OCPs.
\rev{Specifically, for varying initial conditions and horizons, the optimal solutions stay close }to the turnpike during the middle part of the time horizon, whereby with increasing horizon the time spent close to the turnpike grows.
In the well-understood stationary setting, the turnpike corresponds to the optimal steady state.
We refer to \cite{tudo:faulwasser22a} for a recent overview of continuous-time and discrete-time turnpike results.
Indeed, the turnpike phenomenon is closely related to a dissipativity property of the underlying OCP~\cite{Gruene16a,Faulwasser_Automatica17}.
While in the stationary setting the turnpike can easily be precomputed as an optimal steady state, in the time-varying setting the turnpike can turn out to be a more complex trajectory, e.g. an optimal periodic orbit~\cite{Zanon17a}.
In these latter time-varying situations, it is in general very difficult to characterize and precompute the time-varying turnpike without actually solving the OCP in question for a large number of initial conditions. 

In this paper, we provide an analytical characterization of the properties of the time-varying turnpikes that arise in the optimal control formulation of DHN operations.
In particular, we exploit a link between singular arcs arising in the Pontryagin Maximum Principle (PMP) and the exactness of the turnpike property.
The latter means that for varying initial conditions and a sufficiently long horizon, the turnpike is entered exactly.
In~\cite{Faulwasser_EJC17} it was shown that stationary exact turnpikes are closely related to the singular arcs of the OCP.

The main contributions of the present paper are three-fold: 
First, considering an OCP modeling optimal operation of DHNs with time-varying prices and demands, we derive sufficient and necessary conditions for \rev{this} OCP to exhibit a unique time-varying singular arc and we give the closed-form expression.
Second, we \rev{utilize the inherent stability of the DHN model to} extend results on the relation between dissipativity of OCPs and the turnpike property from the stationary setting \cite{Faulwasser_Automatica17} to the time-varying case.
That is, we show that a time-varying exact turnpike implies a strict dissipation \rev{inequality to hold along optimal solutions}, and we show that strict dissipativity of the OCP combined with reachability and a specific structure of the singular arc implies a time-varying exact turnpike.
To the best of our knowledge, our results are the very first to present a structured approach towards the analytic characterization of time-varying turnpike objects in OCPs.
Third, our findings are illustrated with a numerical example.

\subsubsection*{Notation}
For a finite set $\mc{S}$, $|\mc{S}|$ denotes its cardinality.
Let $\R$, $\N$ and $\mbb{C}$ denote the sets of real numbers, positive integers, and complex numbers, respectively.
Moreover, ${\R_{\bullet c}:=\{x\in\R:x\bullet c\}}$ and ${N_{\bullet c}:=\{x\in\N:x\bullet c\}}$ for ${\bullet\in\{<,>,\leq,\geq\}}$ and $c \in \mathbb R$.
For ${x\in\R^n}$, $(x)_i$ denotes its $i^\mrm{th}$ element 
%for ${i\in\{1,\ldots,n\}}$ 
and $\diag{x}$ denotes a diagonal matrix with $x$ on the main diagonal.
For $x_i\in\R$ with $i\in\mc{S}$, ${[x_i]_{i\in\mc{S}}\in\R^{|\mc{S}|}}$ denotes a vector that collects all $x_i$.
The Lebesgue measure on the real line of \rev{a measurable set $\mc{I}$ is denoted by $\mu(\mc{I})$}. 
The set of Lebesgue-integrable functions on $[a,b]\subset\R$ with values in some set $\mc{R}\subseteq\R^{d}$ is $L^1([a,b],\mc{R})$.
Moreover, $W^{1,1}([a,b],\mc{R})$ is the space of functions $x : [a,b]\rightarrow\mc{R}$, such that $x$ and its weak derivative $\dot{x}$ belong to $L^1([a,b],\mc{R})$.
For trajectories ${x\in W^{1,1}([a,b],\mc{R}_1)}$, ${u\in L^1([a,b],\mc{R}_2)}$, we denote the pair of trajectories by ${z=(x,u)\in Z^{1,1}([a,b],\mc{R}_1\times\mc{R}_2)}$.
The class $\mathcal{K}$ denotes functions $\alpha:\R_{\geq0} \rightarrow \R_{\geq0}$ that increase monotonically and are continuous with $\alpha(0)=0$.

%--------------------------------------------------------------------------------------------------------------
\section{Preliminaries and Problem Statement}
%--------------------------------------------------------------------------------------------------------------
We start by summarizing the main aspects of the modeling of DHNs leading to LTI system dynamics.
To this end, we consider DHNs comprising producers and consumers that are interconnected by a closed hydraulic network, i.e., pipes interconnected by junctions.
\subsection{Optimal control formulation for DHNs} 
The following standard assumptions are frequently considered for optimal control of a DHN \cite{Rose_IFAC23,krug_nonlinear_2021,Sibeijn_ECC24,Machado_Automatica22}.

\textit{Modeling Assumptions:}
    (i)~A fixed control volume is associated with each component that is completely filled with water at all times.
    (ii)~The water is incompressible.
    (iii)~Gravitational forces are neglected.
    (iv)~All mass flows are one-dimensional.
    (v)~The internal energy of the water is the main energy source depending linearly on the temperature of the water. Other energy forms are negligible.
    (vi) The mass flows in the DHN are constant.\hfill $\square$

To represent the topology of a DHN, we follow common practice and use a connected graph ${\mc{G}=(\mc{V},\mc{E})}$, where ${\mc{V}=\{v_1,\ldots,v_n\}}$, $n\in\N_{\geq2}$, and ${\mc{E}=\{e_1,\ldots,e_l\}}$, $l\in\N_{\geq2}$, denote the set of vertices and the set of edges, respectively.
Each $v\in\mc{V}$ is assigned to a component. Therefore, 
%with slight abuse of notation, 
we assign to each $v\in\mc{V}$ a temperature $T_v(t)\in\R$, a mass $m_v\in\R_{>0}$, and a heat loss coefficient $\kappa_v\in\R_{>0}$.
Additionally, if $v\in\mc{V}$ is assigned to a producer or a consumer, the heat power $P_v(t)\in\R$ can be injected into or extracted from $v$, respectively. 
To simplify the notation, we set $m_v=1$ for $v\in\mc{V}$, \rev{the ambient temperature}  $T_\mrm{a}=0$ and the heat capacity $c_\mrm{p}=1$ which does not influence the results, i.e., the results are valid for any $m_v>0$, $c_\mrm{p}>0$ \rev{and possibly time-varying $T_\mrm{a}:[0,\infty)\rightarrow\R$}. 
If two components $u,v\in\mc{V}$ are directly hydraulically connected, the edge $e=(u,v)\in\mc{E}$ exists for $u \neq v$.
For an intuitive physical interpretation, we choose the orientation of ${e \in \mc{E}}$ corresponding to the direction of flow of water associated with this edge.
This allows us to define the mass flow on ${e=(u,v)\in\mc{E}}$ as ${q_e\in\R_{>0}}$.
We define the state vector collecting all temperatures as ${x(t)=[T_v(t)]_{v\in\mc{V}}}$, the control input vector collecting injected heat flows as ${u(t)=[P_v(t)]_{v\in\mc{P}}}$, where $\mc{P}\subset\mc V,$ 
${|\mc{P}|=m}$, and the disturbance vector collecting extracted heat flows ${d(t)=[P_v(t)]_{v\in\mc{D}}}$, where $\mc{D}\subset\mc V\setminus\mc P$ and ${w=|\mc{D}|}$, with $m+w\leq n$.
Let ${A_L\in\R^{n\times n}}$ denote the \emph{weighted flow Laplacian} of $\mc{G}$\footnote{For the sake of brevity, we omit a formal definition of the flow Laplacian and refer to \cite{VANDERSCHAFT201721} for details.} using $q_e$ as weight of ${e\in\mc{E}}$, ${A=-A_L-\diag{[\kappa_v]_{v\in\mc{V}}}}$ and ${B\in\R^{n\times m}}$ defined element-wise as
\begin{equation} \label{eq: B}
    (B)_{i,j}=\begin{cases}
        1,\text{ if }v_i\in\mc{V}\text{ is }j^\mrm{th}\text{ producer,}\\
        0,\text{ else.}
    \end{cases}
\end{equation}
Then, the DHN dynamics can be represented by %the LTI system
\begin{equation} \label{eq: LTI system}
    \dot{x}(t) = Ax(t) + Bu(t) + Ed(t),~x(0)=x_0\in\mc{X}_0,
\end{equation}
%models a DHN, 
where ${\mc{X}_0\subset\R^n}$ denotes a compact set of initial conditions, and ${E\in\R^{n\times w}}$ analogously to \eqref{eq: B} for $w$ consumers.

Moreover, we consider input constraints given by upper limits and lower limits for $u$, i.e., \rev{${\mc{U}=[(u_\mrm{min})_1,(u_\mrm{max})_1]\times\cdots\times[(u_\mrm{min})_m,(u_\mrm{max})_m]\subset\R_{\geq0}^m}$.}
The considered time-varying operating costs are
\begin{equation} \label{eq: ell}
    \ell(\tau,x,u) = {\textstyle \frac{1}{2}} x^\top Q x + u^\top S x + x^\top r(\tau) + u^\top p(\tau),
\end{equation}
where $Q\in\R^{n\times n}$ is a symmetric matrix, $S\in\R^{n\times m}$, ${r:[0,\infty)\rightarrow\mathbb{R}^n}$ and ${p:[0,\infty)\rightarrow\mathbb{R}^m}$.

With these considerations, we arrive at the following OCP. % modeling optimal operation of DHNs:
Given $T\in(0,\infty)$, ${x(0)=x_0\in\mc{X}_0}$, \rev{and functions $r, p, d$,} find the control input trajectory $u_T^\star$ that solves
\begin{align*}\label{prob: OCP}%J_T(x_0,u) =
    \min_{u\in L^1([0,T],\mc{U})}&  \int_{0}^{T} \ell(\tau,x(\tau),u(\tau))\, {\rm d}\tau \tag{$\mc{OCP}_T(x_0)$}\\
    \text{s.t.}~&\text{\eqref{eq: LTI system}},~u(\tau)\in\mc{U},
\end{align*}
for all $\tau\in[0,T]$.\hfill $\square$\\
Observe that in \ref{prob: OCP} the time horizon $T$ is fixed in the optimization, but we are interested in characterizing the solution for different $T\in(0,\infty)$ and different initial conditions $x_0 \in \mathcal X_0$. As a shorthand highlighting the parametric nature of the problem, we use \ref{prob: OCP}.
The optimal state trajectory resulting from the optimal input $u_T^\star(\cdot, x_0)$ is denoted by $x_T^\star(\cdot;x_0)$.
The optimal pair is written as ${z_T^\star(\cdot;x_0)=(x_T^\star(\cdot;x_0),u_T^\star(\cdot;x_0))}$ and we have the set ${\mc{Z}=\R^n \times \mc{U}}$.
\rev{To simplify notation, the dependency on $x_0$ is dropped when not explicitly needed.}

\begin{assumption} \label{ass: technical modeling}
    The functions $d$ from \eqref{eq: LTI system} and $r$, $p$ from \eqref{eq: ell} are known and bounded for all $t\geq0$ and $k$-times continuously differentiable for $k>0$ sufficiently large.\footnote{The minimum required smoothness can be determined explicitly by considering the so-called DAE index of the system \eqref{eq: NCO}~\cite{KunkelMehrmann_DAE}. Details are omitted for the sake of brevity.} 
    \hfill $\square$
\end{assumption}
Assumption~\ref{ass: technical modeling} implies that system \eqref{eq: LTI system} has a unique solution $x\in W^{1,1}([0,T],\R^n)$ for every $u\in L^1([0,T],\R^m)$.

\subsection{Recap -- Singular arcs in OCPs}
The Hamiltonian of \ref{prob: OCP} reads
\begin{equation}
    H(\tau,\lambda_0,\lambda,x,u) = \begin{bmatrix}
        \lambda_0 \\ \lambda
    \end{bmatrix}^\top \begin{bmatrix}
        \ell(\tau,x,u) \\ Ax + Bu + Ed(\tau)
    \end{bmatrix},
\end{equation}
where ${\lambda_0\in\R}$ and ${\lambda\in\R^n}$ are the costates/adjoints for cost and dynamics.
The PMP gives that for any ${z_T^\star(\cdot,x_0)}$ solving \ref{prob: OCP}, there is $\lambda\in W^{1,1}([0,T],\R^n)$ such that
\begin{subequations} \label{eq: PMP}
    \begin{align}
        \dot{x}_T^\star(\tau) &= Ax_T^\star(\tau) + Bu_T^\star(\tau) + Ed(\tau),~x_T^\star(0)=x_0, \label{eq: PMP state}\\
        \dot{\lambda}(\tau) &= -(Qx_T^\star(\tau)+S^\top u_T^\star(\tau)+r(\tau))%-\frac{1}{2}Q\bar x) 
         \rev{- A^\top\lambda(\tau),} \label{eq: PMP costate}\\&~~~~\rev{\lambda(T)=0},\nonumber  \\
        u_T^\star(\tau) &= \underset{u\in\mc{U}}{\mrm{argmin}}~u^\top (S x_T^\star(\tau) + B^\top \lambda(\tau) + p(\tau)) \label{eq: stationary condition}.
    \end{align}
\end{subequations}
The \ref{prob: OCP} does not involve any terminal constraints.
Thus, ${\lambda(T) = 0}$ for ${T< \infty}$ and the PMP only holds if ${\lambda_0 \not = 0}$~\cite{Carlson91}.
Hence, we set $\lambda_0=1$ without loss of generality.\footnote{In other words, for all $T< \infty$, the considered OCP is not abnormal. Moreover, we remark that for $T= \infty$ the transversality condition of the costate has to be dropped, see \cite{tudo:faulwasser21a} for recent results how to characterize $\lambda(T=\infty)$ in a dissipativity framework.}

Since \eqref{eq: stationary condition} is affine linear in $u$, i.e., the OCP is singular, we consider the switching function of the $i$-th input
\begin{equation} \label{eq: switchFunc}
    s_i(\tau) = \big(S x_T^\star(\tau) + B^\top \lambda(\tau) + p(\tau)\big)_i
\end{equation}
for all $i=1,\ldots,m$ and all ${\tau\in[0,T]}$.
Consider the set \rev{${\mc{T}_{i,\bullet}(T) =\{\tau\in[0,T]:s_i(\tau)\bullet0\}}$} with ${\bullet\in\{>,<,=\}}$. 
It follows that ${(u_T^\star(\tau))_i=(u_\mrm{min})_i}$ for ${\tau \in \mc{T}_{i,>}(T)}$ and ${(u_T^\star(\tau))_i=(u_\mrm{max})_i}$ for ${\tau \in \mc{T}_{i,<}(T)}$\rev{, i=1,\ldots,m}. %, where $\mu(\cdot)$ is the Lebesgue measure.
\rev{Based on \cite{Faulwasser_EJC17}, we slightly adapt the definition of \emph{singular arcs} to the time-varying setting considered in this work.}
\begin{defi}[Unique singular arc] \label{def: singular arc}
    \ref{prob: OCP} is said to exhibit a \emph{unique singular arc} if \rev{$\mu\left(\bigcap_{i=1}^m \mc{T}_{i,=}(T)\right)>0$, and the tuple $(z_T^\star(\tau),\lambda(\tau))$ solving \eqref{eq: PMP} on $\tau\in\bigcap_{i=1}^m \mc{T}_{i,=}(T)$} is unique almost everywhere. \hfill $\square$
\end{defi}
In other words, Definition \ref{def: singular arc} demands that the solution to \ref{prob: OCP} is unique almost everywhere when no input constraint is active. %, i.e., $u_T^\star(\cdot)\notin\mrm{int}(\mc{U})$.
Hence, to analyze whether \ref{prob: OCP} exhibits a unique singular arc\rev{---simultaneously for all control inputs $u_1,\dots, u_m$---}in the sense of Definition \ref{def: singular arc},\footnote{\rev{The extension to the case of only a subset of the inputs $u_1,\dots, u_m$ taking on values in the interior of the constraint set $\mathcal U$ is beyond the scope of the present work.}} we examine the optimality system \eqref{eq: PMP} on singular arcs, i.e., when no input constraint is active.
To this end, we define ${\xi^\star(\tau)=\begin{bmatrix}{x_T^\star(\tau)}^\top & \lambda(\tau)^\top & {u_T^\star}(\tau)^\top\end{bmatrix}^\top}$, ${f(\tau)=\begin{bmatrix}{(Ed(\tau))}^\top & -r(\tau)^\top & {p(\tau)}^\top\end{bmatrix}^\top}$.
When no input constraint is active, \eqref{eq: PMP} leads to the DAE
\begin{equation} \label{eq: NCO}
    D\dot{\xi^\star}(\tau)=M\xi^\star(\tau)+f(\tau),
\end{equation}
where
\begin{equation*}
    \begin{aligned}
        D&=\begin{bmatrix}
            I_n & 0 & 0 \\ 0 & I_n & 0 \\ 0 & 0 & 0
        \end{bmatrix},~
        M=\begin{bmatrix}
            \phantom{-}A\phantom{^\top}& \phantom{-}0\phantom{^\top} & \phantom{-}B\phantom{^\top} \\ -Q\phantom{^\top} & -A^\top & -S^\top \\ \phantom{-}S\phantom{^\top} & \phantom{-}B^\top & \phantom{-}0\phantom{^\top}
        \end{bmatrix},
    \end{aligned}
\end{equation*}
and where $0$ represents the zero matrix with the appropriate sizes. 
Note that the initial condition of this DAE is \rev{implicitly} specified by the switching functions of all inputs, \eqref{eq: switchFunc}, being zero.
We recall the following standard assumptions for constant-coefficient DAEs \cite{KunkelMehrmann_DAE}.

\begin{assumption} \label{ass: regular pencil}
The matrix pencil $sD-M, s\in \mathbb C$ is \emph{regular}, i.e., $\det(sD-M)$ is not the zero polynomial. \hfill $\square$
\end{assumption}
To analyze the existence and uniqueness of solutions to \eqref{eq: NCO}, we recall standard results for linear DAEs with constant coefficients.  % that are used in what follows.

\begin{proposition}[Lem. 2.8 and Thm. 2.12 of \cite{KunkelMehrmann_DAE}] \label{prop: DAE standard}
    Let Assumption~\ref{ass: regular pencil} hold.
    Then, \eqref{eq: NCO} has a unique solution for any consistent initial value and can be transformed into \emph{Weierstraß canonical form}.
    That is, there exist nonsingular matrices $P_\mrm{w}$, $Q_\mrm{w}$, a nilpotent matrix $N_\mrm{w}$, and a matrix $J_\mrm{w}$ in \emph{Jordan canonical form}, all in appropriate dimensions, such that the unique solution of \eqref{eq: NCO} is given by
    %\begin{subequations}
    \begin{align}
    \label{eq: singArcSol}
        \Tilde{\xi}_2(\tau)= -\!\!\!\!\!\sum_{i=0}^{2n+m-1}\!\!\! N_\mrm{w}^i \Tilde{f}_2^{(i)}(\tau),~ 
        \dot{\Tilde{\xi}}_1(\tau)= J_\mrm{w} \Tilde{\xi}_1(\tau) + \Tilde{f}_1(\tau),
    \end{align}
    where $\begin{bmatrix} \Tilde{\xi}_1(\tau) \\ \Tilde{\xi}_2(\tau) \end{bmatrix} = Q_\mrm{w}^{-1} \xi^\star(\tau),~\begin{bmatrix} \Tilde{f}_1(\tau) \\ \Tilde{f}_2(\tau) \end{bmatrix} = P_\mrm{w} f(\tau)$.\hspace*{1.3cm} $\square$ 
\end{proposition}

\subsection{Dissipativity of OCPs}
\rev{Let $\bar z$ denote an arbitrary infinite-horizon input-state pair satisfying (2) and the input constraints and c}onsider the \rev{shifted} supply rate ${\ell_{\bar z}(\tau,z)=\ell(\tau,z)-\ell(\tau,\bar z(\tau))}$.
\begin{defi}[Time-varying strict dissipativity of OCPs] \label{def: strict dissipativity}
    \ref{prob: OCP} is said to be \emph{strictly dissipative} with respect to ${\bar z\in \rev{W}^{1,1}(\rev{[0,\infty)},\mc{Z})}$ if there exists a nonnegative storage function $\mc{S}:\R\times\R^n\rightarrow\R_{\geq s}$ for some $s\geq 0$ such that, for all $x_0\in\mc{X}_0$, $T\geq0$ and along optimal pairs $z_T^\star(\cdot, x_0)$ the inequality
    \begin{multline} \label{eq: dissipation inequality}
        \mc{S}(T,x_T^\star(T))-\mc{S}(0,x_0) 
        \leq \int_{0}^{T}\hspace*{-0.25cm} -\alpha\big(\|z_T^\star(\tau, x_0)-\bar z(\tau)\|\big)\\+\ell_{\bar z}(\tau, z_T^\star(\tau, x_0)) ~\rm \rm d\tau,
        \tag{sDI}
    \end{multline}
    holds for some ${\alpha\in\mc{K}}$.\hfill $\square$
\end{defi}

Observe that, similar to \cite{kit:faulwasser18c,Faulwasser_Automatica17}, this notion of dissipativity requires the strict dissipation inequality \eqref{eq: dissipation inequality} to hold only along optimal solutions. Jan C. Willems provided a dissipativity characterization via the available storage~\cite{Willems72a}.
The available storage is adapted to strict dissipativity by introducing $\ell_{{\bar z},\alpha}(\tau,z):=\ell_{\bar z}(\tau,z)-\alpha\big(\|z-\bar z(\tau)\|\big)$ and given as the value function of the free end-time OCP
\begin{align}\label{prob: Verify dissipativity}
        \mc{S}^\mrm{a}(x_0)= &\sup_{u\rev{_T^\star}\in L^1([0,T],\mc{U}),~T\geq 0} -\int_{0}^{T} \ell_{{\bar z},\alpha}(\tau,z(\tau))~\rm d\tau  \\
        &\text{s.t.}~\text{\rev{$u_T^\star$ solves \ref{prob: OCP}},} %\eqref{eq: LTI system}},~x(0)=x_0,~%u(\tau)\in\mc{U}},
        %~T\geq0, 
        \nonumber
    \end{align}
see also \cite{Faulwasser_Automatica17,Willems72a,cherifi2024}.  
The following result recalls the available storage characterization of strict dissipativity~\cite{Willems72a}.

\begin{theorem} \label{thm: Willems}
    \ref{prob: OCP} is strictly dissipative if and only if $\mc{S}^\mrm{a}(x_0)$ from~\eqref{prob: Verify dissipativity} is finite for all $x_0 \in \mc{X}_0$.\hfill $\square$
\end{theorem}

%--------------------------------------------------------------------------------------------------------------
\section{Main Results}
%--------------------------------------------------------------------------------------------------------------
The following result is a direct consequence of Proposition~\ref{prop: DAE standard} on the existence of a unique solution to linear DAEs with constant coefficients. Henceforth, it will be used to characterize exact turnpikes in \ref{prob: OCP}.
\begin{corollary}[Unique singular arcs in \ref{prob: OCP}] \label{cor: unique singular arc}
    Let Assumption \ref{ass: technical modeling} hold and let $\tau\in\bigcap_{i=1}^m \mc{T}_{i,=}(T)\neq\emptyset$.
    Then, \ref{prob: OCP} exhibits a unique singular arc in the sense of Definition \ref{def: singular arc} if and only if Assumption \ref{ass: regular pencil} holds. Moreover, on the singular arc, the solution is given by \eqref{eq: singArcSol}.\hfill $\square$
\end{corollary}

%--------------------------------------------------------------------------------------------------------------
\subsection{Time-varying exact turnpikes for DHNs}
%--------------------------------------------------------------------------------------------------------------
Next, we extend the \rev{dissipativity analysis  for stationary turnpikes in \ref{prob: OCP} from \cite{Faulwasser_Automatica17}} to the time-varying case. We provide sufficient conditions for \ref{prob: OCP} to exhibit a time-varying exact turnpike property and we show that a time-varying exact turnpike implies strict dissipativity.% in the sense of Definition~\ref{def: strict dissipativity}.

First, following \cite{tudo:faulwasser22a,Faulwasser_Automatica17}, we recall the definition of a time-varying measure turnpike property. 
Consider an arbitrary infinite-horizon input-state pair $\bar z$  satisfying \eqref{eq: LTI system} and the input constraints and the set %\hg{HG: Die Menge hängt noch vom optimal pair ab, aber das muss doch nicht eindeutig sein? Wie geht man damit um?}
    \begin{equation} \label{eq: defi Theta}
    \Theta_T(\varepsilon):=\{\tau\in[0,T]:\|z_T^\star(\tau, x_0)-\bar{z}(\tau)\|>\varepsilon\},
\end{equation}
    This set is Borel measurable and captures all time points when an optimal pair $z_T^\star$ is outside of an $\varepsilon$ neighborhood of $\bar z$.
\begin{defi}[\rev{Time-varying measure turnpike}] \label{def: turnpike}
     \ref{prob: OCP} exhibits a time-varying turnpike property with respect to ${\bar z}$ if there exists a function ${\nu:[0,\infty)\rightarrow[0,\infty)}$
     such that, for all ${x_0\in\mc{X}}$, all $T>0$, \rev{all optimal $u_T^\star\in L^1([0,T],\mathbb{R}^m)$}, and all $\varepsilon> 0$ we have
    \begin{equation} \label{eq: turnpike}
\mu(\Theta_T(\varepsilon))\leq\nu(\varepsilon)<\infty. 
    \end{equation}
    The optimal pairs $z_T^\star(\cdot, x_0)$ exhibit a time-varying \textit{exact turnpike property} if \eqref{eq: turnpike} holds for $\varepsilon=0$.\hfill $\square$
\end{defi}
The Lebesgue measure $\mu(\Theta_T(\varepsilon))$ quantifies the size $\Theta_T(\varepsilon)$, i.e., the time spent off the time-varying turnpike.
Put differently, Definition \ref{def: turnpike} requires that the time $\mu(\Theta_T(\varepsilon))$  the optimal solutions to \ref{prob: OCP} spend outside an $\varepsilon$-neighborhood of $\bar z$ is bounded by $\nu(\varepsilon)$, where $\nu(\varepsilon)$ is independent of the horizon $T$.
This essentially means that for any $x_0 \in \mc{X}_0$ and for sufficiently large $T>0$, the optimal trajectory $z_T^\star(\cdot, x_0)$ has to enter any arbitrarily small $\varepsilon$-neighborhood of $\bar z$. In the case of $\varepsilon=0$, i.e., the time-varying exact turnpike property, $z_T^\star(\cdot, x_0)$ must be temporarily equal to the turnpike.

Similarly to the results for time-invariant OCPs exhibiting stationary turnpikes, a time-varying turnpike of \ref{prob: OCP} can be shown if strict dissipativity is given.
Hence, we have the following assumption.
\begin{assumption} \label{ass: reachability and dissipativity}
(i) \ref{prob: OCP} is strictly dissipative with respect to ${\bar z\in Z^{1,1}(\rev{[0,\infty)},\mrm{int}(\mc{Z}))}$ in the sense of Definition~\ref{def: strict dissipativity}.
    (ii)  The trajectory ${\bar z\in Z^{1,1}([0,\infty),\mrm{int}(\mc{Z}))}$ is exponentially reachable, i.e., for all $x_0\in\mc{X}_0$, there exist
    an infinite-horizon admissible input ${u_\infty\in L^1([0,\infty),\mc{U})}$, and \rev{${c_1,~c_2\in\R_{>0}}$}, independent of $x_0$, such that, \rev{for all $\tau\geq 0$}, ${\|(x(\tau;u_\infty),u_\infty(\tau))-\bar z(\tau)\|\leq c_1 e^{-c_2\tau}}$ holds.
\end{assumption}

Next, we extend \cite[Thm. 2]{Faulwasser_Automatica17} to the time-varying case.

\begin{proposition}[Strict dissipativity $\Rightarrow$ turnpike] \label{prop: dissipativity to turnpike} 
Let Assumption \ref{ass: reachability and dissipativity} hold.
    Then, \ref{prob: OCP} exhibits a time-varying turnpike property with respect to ${\bar z\in Z^{1,1}([0,T],\mrm{int}(\mc{Z}))}$.\hfill $\square$
\end{proposition}

The proof follows the same steps as the one of \cite[Thm. 2]{Faulwasser_Automatica17} and is omitted due to space limitations.
The next result is an extension of \cite[Thm. 3]{Faulwasser_EJC17} to time-varying exact turnpikes.
\begin{theorem}[Turnpikes in \ref{prob: OCP} are exact]
    Let Assumptions \ref{ass: technical modeling}, \ref{ass: regular pencil} and \ref{ass: reachability and dissipativity} hold. 
    Then, for all $x_0\in\mc{X}_0$ and sufficiently large $T>0$, the turnpike of \ref{prob: OCP} 
    %has to be 
    is exact in the sense of Definition \ref{def: turnpike}.\hfill $\square$
\end{theorem}
\begin{proof}
    Let $\mc{B}_\varepsilon(\bar{z}(\tau))$ denote a closed ball centered at the turnpike at time $\tau$ with radius $\varepsilon>0$.
    From Proposition~\ref{prop: dissipativity to turnpike}, $\mu(\Theta_T(\varepsilon))\leq\nu(\varepsilon)$ follows.
    Thus, ${z_T^\star(\tau)\in\mc{B}_\varepsilon(\bar{z}(\tau))}$ \rev{for all}  ${\tau\in[0,T]\setminus\Theta_T(\varepsilon)}$ and \rev{all} $T>0$. % \todo{R1: "sufficiently large $T > 0$"; however, either it is not clear why this is mentioned, or the argument demands further clarifications.}.
    Let ${\mc{B}_\varrho(\bar{z}(\tau))}$ denote another closed ball centered at $\bar{z}$ at time $\tau$ with radius $\varrho>0$, such that ${\mc{B}_\varrho(\bar{z}(\tau))\subset\mrm{int}(\mc{Z})}$ for ${\tau\in[0,T]}$.
    We define %\todo{R1: The role of $\hat{\epsilon}$ is not clear.}
    \begin{equation} \label{eq: hat epsilon}
        \hat{\varepsilon}=\min_{\tau\in[0,T]}\{\|u_\mrm{min}-\bar{u}(\tau)\|,\|u_\mrm{max}-\bar{u}(\tau)\|,\varrho\}.
    \end{equation}
    \rev{Assumption \ref{ass: reachability and dissipativity} shows that close to the turnpike (in terms of $z^\star_T$) we have $\bar u(\tau)\in\mrm{int}(\mc{U})$} (almost surely) and hence ${\hat{\varepsilon}>0}$ holds. 
    \rev{Consequently, Proposition~\ref{prop: dissipativity to turnpike} yields that} for some ${\varepsilon\in(0,\hat{\varepsilon})}$, sufficiently large ${T>0}$, and ${\tau\in[0,T]\setminus\Theta_T(\varepsilon)}$, we have ${z_T^\star(\tau)\in\mc{B}_\varepsilon(\bar{z}(\tau))}$. Thus, we henceforth analyze \ref{prob: OCP} locally around $\bar{z}$.
    For $\varepsilon\in[0,\hat{\varepsilon})$, $\tau\in[0,T]\setminus\Theta_T(\varepsilon)$, we show  $z_T^\star(\tau)\in\mc{B}_\varepsilon(\bar{z}(\tau))$ implies $z_T^\star(\tau)=\bar{z}(\tau)$ (Steps 1 and 2).
    Step~3 shows $\mu(\Theta_T(0))\leq\nu(\hat{\varepsilon})<\infty$.
    
    Step 1: We show that for $\varepsilon\in[0,\hat{\varepsilon})$ and $\tau\in[0,T]\setminus\Theta_T(\varepsilon)$
      \begin{equation} \label{eq: step 1}
        z_T^\star(\tau)\in\mc{B}_\varepsilon(\bar z(\tau))\quad \Longrightarrow\quad  u_T^\star(\tau)=\bar u(\tau).
    \end{equation}
    
    Since $\mc{B}_\varepsilon(\bar z(\tau))\subset\mc{B}_\varepsilon(\bar x(\tau))\times\mc{B}_\varepsilon(\bar u(\tau))$ for $\varepsilon\in[0,\hat{\varepsilon})$, \rev{we have that $u_T^\star(\tau)\in \mc{B}_\varepsilon(\bar u(\tau))$ holds for all $\varepsilon\in[0,\hat{\varepsilon})$ and all $\tau\in[0,T]\setminus\Theta_T(\varepsilon)$. } 
    For $u_T^\star$ solving \ref{prob: OCP}, \rev{and for $\tau\in[0,T]\setminus \mathcal{T}_{i,=}(T)$ for some $i=1,\ldots,m$, we obtain from the PMP that $u_T^\star\in\{u_{\min,i},u_{\max,i}\}$ holds.
    Contrary to that, suppose that $\tau\in\bigcap_{i=1}^m\mathcal{T}_{i,=}(T)$, then the uniqueness of the singular arc (Corollary~\ref{cor: unique singular arc}) implies that  $u_T^\star(\tau)=\overline{u}(\tau)$ holds.
    Consequently, for all $\tau\in[0,T]$}, ${(u_T^\star(\tau))_i \in \{u_{\mrm{min},i}, u_{\mrm{max},i}, \bar u_i(\tau)\},~\text{for all } i=1,\ldots,m,}$ holds.
    That is, the optimal control is component-wise either at its max or min value or at the turnpike. 
    Hence, for all ${\varepsilon\in[0,\hat{\varepsilon})}$ and ${\tau\in[0,T]\setminus\Theta_T(\varepsilon)}$, \eqref{eq: hat epsilon} implies
    \begin{multline*}
        \mc{B}_\varepsilon(\bar u(\tau)) \cap \big( \{u_{\mrm{min},1}, u_{\mrm{max},1}, \bar u_1(\tau)\} \times \cdots \\
        \times \{u_{\mrm{min},m}, u_{\mrm{max},m}, \bar u_m(\tau)\}\big) = \bar u(\tau).
    \end{multline*}
    Put differently, strict dissipativity of the OCP combined with the uniqueness of the singular arc gives that sufficiently close to the turnpike the input can only be at its turnpike value \rev{$\bar u(\tau)$}. Therefore, \eqref{eq: step 1} holds.

    Step 2: From ${u_T^\star(\tau) = \bar u(\tau) \in \mrm{int}(\mc{U})}$ for ${\tau\in[0,T]\setminus\Theta_T(\varepsilon)}$, we conclude that \rev{$\mu\left(\bigcap_{i=1}^m \mc{T}_{i,=}(T)\right)>0$}. % \todo{R1: conclusion not clear}.
    Due to Assumption~\ref{ass: regular pencil}, \ref{prob: OCP} exhibits a unique singular arc for ${\tau\in[0,T]\setminus\Theta_T(\varepsilon)}$. Hence, it follows that
    \begin{equation} \label{eq: step 2}
        z_T^\star(\tau)\in\mc{B}_\varepsilon(\bar z(\tau))\Rightarrow x_T^\star(\tau) = \bar x(\tau).
    \end{equation}
    Combining \eqref{eq: step 1} and \eqref{eq: step 2} gives
    \begin{equation} \label{eq: z at turnpike}
        z_T^\star(\tau)\in\mc{B}_\varepsilon(\bar{z}(\tau)) \Rightarrow z_T^\star(\tau)=\bar{z}(\tau).
    \end{equation}
    Step 3: Proposition~\ref{prop: dissipativity to turnpike} implies $\mu(\Theta_T(\varepsilon))\leq\nu(\varepsilon)<\infty$ for all $\varepsilon>0$. 
    From \eqref{eq: z at turnpike}, we can deduce that, for all ${\varepsilon_1, \varepsilon_2 \in [0,\hat{\varepsilon})}$, it holds that \rev{$\mu(\Theta_T(\varepsilon_1))=\mu(\Theta_T(\varepsilon_2))<\nu(\hat{\varepsilon})$} 
    that is, $\mu(\Theta_T(0))\leq\nu(\hat{\varepsilon})<\infty$.
\end{proof}

\subsection{A converse turnpike result}
A crucial breakthrough in the dissipativity-based turnpike analysis have been converse results, i.e., results which explore the implications of assuming the existence of a turnpike in an OCP, see, e.g., \cite{Gruene16a,Faulwasser_Automatica17}.
The next result, which is a time-varying generalization of \cite[Thm. 7]{Faulwasser_Automatica17} shows that exact time-varying turnpikes imply strict dissipativity.
\begin{proposition}[Exact turnpike $\Rightarrow$ strict dissipativity] \label{prop: turnpike to dissipativity}
    \rev{L}et \ref{prob: OCP} exhibit a time-varying exact turnpike in the sense of Definition \ref{def: turnpike} at ${\bar z\in Z^{1,1}(\rev{[0,\infty)},\mrm{int}(\mc{Z}))}$. %Let Assumption~\ref{ass: exact turnpike} hold.
    
    \rev{If $\bar z$ is bounded, }
    then \ref{prob: OCP} is strictly dissipative with respect to ${\bar z}$ \rev{and all $\alpha\in\mathcal{K}$}.\hfill $\square$
\end{proposition}
\begin{proof}
    Choose any $\alpha\in\mc{K}$, $x_0\in\mathcal{X}_0$ and an optimal control $u_T^\star(x_0)\in L^1([0,T],\mc{U})$ for \ref{prob: OCP}. Let $\varepsilon=0$ and split the horizon $[0,T]$ into $\Theta_T(0)$ and $\Delta = [0,T]\setminus\Theta_T(0)$.
    We aim to derive an upper bound on the integral in the available storage $\mc{S}^\mrm{a}(x_0)$ that is given by \eqref{prob: Verify dissipativity}:
    \begin{align}
    &~~~~-\int_0^T\ell_{\bar z,\alpha}(\tau,z_T^\star(\tau))\mathrm{d}\tau \nonumber \\&
    =-\int_{\Theta_T(0)}\ell_{\bar z,\alpha}(\tau,z_T^\star(\tau))\mathrm{d}\tau -\int_\Delta \ell_{\bar z,\alpha}(\tau,z(\tau))\mathrm{d}\tau \nonumber \\&\leq\int_{\Theta_T(0)}|\ell_{\bar z,\alpha}(\tau,z_T^\star(\tau))|\mathrm{d}\tau + \int_\Delta |\ell_{\bar z,\alpha}(\tau,z(\tau))|\mathrm{d}\tau \nonumber \\
&    \rev{=} \int_{\Theta_T(0)}|\ell_{\bar z,\alpha}(\tau,z_T^\star(\tau))|\mathrm{d}\tau,\label{eq:first_estimate}
    \end{align}
    where in the last estimate the second integral vanishes as ${z_T^\star = \bar z}$  on ${\Delta = [0,T]\setminus\Theta_T(\varepsilon)}$.
    Using Lemma~\ref{lem: bounded states} we find that $x_T^\star$ admits an upper bound that does not depend on $T$.
    Hence, there exists $0<\hat \ell<\infty,$ such that\vspace*{-2mm}
    \begin{align}
    &~~      |\ell_{\bar z,\alpha}(\tau,z_T^\star(\tau))|\nonumber \\
&\hspace{-.15cm}=|\ell(\tau,z_T^\star(\tau))-\ell(\tau,\bar z(\tau))-\alpha(\|z_T^\star(\tau)-\bar z(\tau)\|)| \leq \hat\ell \label{eq:second_estimate}\vspace*{-1mm}
    \end{align}
is valid for all $\tau\in[0,T]$ and $\hat \ell$ is independent of both $T\geq 0$ and $u_T^\star$. Combining \eqref{eq:first_estimate}, \eqref{eq:second_estimate} and \eqref{eq: turnpike}, we obtain
\begin{align}
    \label{eq:almost}
-\int_0^T\ell_{\bar z,\alpha}(\tau,z_T^\star(\tau))\mathrm{d}\tau\leq \mu(\Theta_T(0))\hat\ell \leq \nu(0)\hat \ell
\end{align}
for all $T\geq 0$ and all solutions  $u_T^\star\in L^1([0,T],\mc{U})$ of \ref{prob: OCP}. Taking now the supremum in \eqref{eq:almost} over all optimal controls $u_T^\star\in L^1([0,T],\mc{U})$ and all $T\geq 0$ does, by construction, not affect the upper bound $\hat \ell$ and leads to ${\mc{S}^\mrm{a}(x_0)\leq \nu(0)\hat \ell<\infty}$ for all $x_0\in\mc{X}_0$.
Using Theorem~\ref{thm: Willems}, the assertion follows. 
\end{proof}
\vspace*{-2mm}

%--------------------------------------------------------------------------------------------------------------
\section{Numerical example} \label{sec: numerical example}
\begin{figure}
    \centering
    \includegraphics[width=\columnwidth]{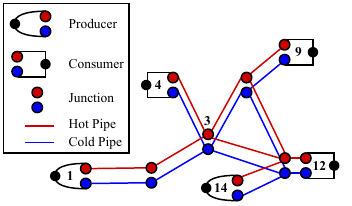}
    \caption{Graph $\mc{G}$ of example DHN with \rev{vertex numbering}.}
    \label{fig: net_plot}
    \vspace*{-5mm}
\end{figure}

%--------------------------------------------------------------------------------------------------------------
\rev{As a numerical example, we consider the DHN shown schematically Figure \ref{fig: net_plot}.
A detailed derivation of the numerical parameters is omitted for brevity\footnote{The source code can be accessed 
 at \url{https://github.com/max65945/Turnpikes-of-DHNs}.} and only parameters} required to verify Assumption \ref{ass: regular pencil} are defined below.

\begin{figure}
    \centering
    \includegraphics[width=\columnwidth]{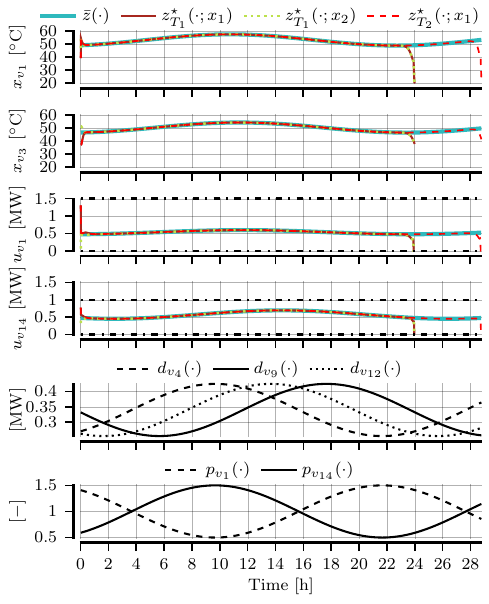}
    \caption{Simulation results of \ref{prob: OCP} for exemplary DHN. % with ${T_1=24\mrm{h}}$, ${T_2=29\mrm{h}}$, ${x_1=0.8x_\mrm{n}(0)}$ and ${x_2=1.1x_\mrm{n}(0)}$.
    The dash dotted line represents $\mc{U}$.}
    \label{fig: sim_plot}
    \vspace*{-5.5mm}
\end{figure}

We select ${S=B^\top}$, ${Q\in\R^{15\times15}}$ as diagonal matrices with all diagonal entries equal to $10^3$, $r(\tau)=-Qx_\mrm{n}$, where $x_\mrm{n}$ denotes nominal temperatures, and sinusoidal $p(\tau)$ and $d(\tau)$ as in Figure~\ref{fig: sim_plot}.
Inserting the defined parameters into \eqref{eq: NCO} yields a nonzero polynomial for $\det(sD-M)$ satisfying Assumption~\ref{ass: regular pencil}. 
With the given parameters, \ref{prob: OCP} is solved \rev{numerically} for the horizons $T_1=24\mrm{h}$ and $T_2=29\mrm{h}$, and for the initial states $x_1=0.8x_\mrm{n}$ and $x_2=1.1x_\mrm{n}$.
The simulation results for the temperature at a producer with the associated heat flow injection and the temperature at a consumer are shown in Figure \ref{fig: sim_plot}.
The turnpike property can be clearly observed for different initial conditions and time horizons.

\section{Conclusions and Outlook}
Motivated by the objective to \rev{contribute to decarbonizing the heat sector through the implementation of MPC} for DHNs, we have conducted an analysis of time-varying turnpikes in \rev{the underlying} OCPs.
In particular, we have analyzed the underlying singularity structure of the OCP to show when arising time-varying turnpikes are exact.
We have also presented a converse turnpike result which shows that, when an exact time-varying turnpike arises, the underlying OCP is strictly dissipative. 
\rev{Future work will analyze the DHN structure in more detail and use the results above for  MPC design.}
Moreover, the extension to settings with varying mass flows appears promising.
%inherent passivity properties

\section*{Appendix}
\begin{lemma}\label{lem:AHurwitz}
    The matrix $A$ from \eqref{eq: LTI system} is Hurwitz.\hfill $\square$
\end{lemma}
\begin{proof}
    %By definition of the graph Laplacian in \cite{VANDERSCHAFT201721}, and 
    Since \rev{$m_v=1$ for all ${v\in\mc{V}}$ and }$\diag{[\kappa_v]_{v\in\mc{V}}}$ %\todo{R6: In the proof of Lemma 1, there is $m_v$ appearing, but $m_v=1$ for all v. This is also why $A$ can be chosen as the graph Laplacian without including weights corresponding to $m_v$. It would be consistent to argue here with the identity matrix I instead of $\diag{[m_v]_v}$.} 
    is diagonal positive definite, we have that $A$ is strictly diagonally dominant with negative diagonal entries.
    Thus, invoking Gershgorin's circle theorem \cite[Theorem 6.2.8]{HorJ12}, we conclude that $A$ is Hurwitz.
\end{proof}
Assumption~\ref{ass: technical modeling} (iii) implies the following, cf. \cite[p. 94]{khalil}:
\begin{lemma}[Bounded states] \label{lem: bounded states}
    For ${x_0\in\mc{X}_0}$ and ${u\in L^1([0,T],\mc{U})}$, and consider \eqref{eq: LTI system} with $A$ Hurwitz.
    Let $\omega(A)<0$ denote an upper bound for the real-parts of the eigenvalues of $A$.
    From Assumption \ref{ass: technical modeling} and the definition of $\mc{U}$ follows that ${\hat d=\sup_{t\in[0,T]} \|d(t)\|<\infty}$ and ${\hat{u}=\sup_{u\in\mc{U}}\|u\|<\infty}$. 
    Then, for some $k>0$ the state trajectory of \eqref{eq: LTI system} is bounded for all $T\geq0$, i.e., ${\| x(T) \|\leq ke^{-|\omega(A)|T}~\|x_0\| + \textstyle{\frac{k}{|\omega(A)|}} (\|B\| \hat u + \|E\| \hat d). \hspace*{.5cm} \square}$
\end{lemma}

\bibliographystyle{IEEEtran}
\bibliography{refs}

\end{document}